\documentclass[12pt,]{article}
\usepackage{amsfonts,amssymb,amsthm,eucal,amsmath,mathtools}
\usepackage{colortbl}
\usepackage{xy}
\input xy
\xyoption{all}
\allowdisplaybreaks
\newtheorem{thm}{Theorem}
\newtheorem{cor}[thm]{Corollary}
\newtheorem{lemma}[thm]{Lemma}

\newtheorem{prop}[thm]{Proposition}

\date{}

\begin{document}

\title{A note on Mahler's conjecture
}

\author{Shlomo Reisner, Carsten Sch\"utt and Elisabeth M. Werner
\thanks{Part of the work was done while the three authors visited the  R´enyi Institute
of Mathematics of  the Hungarian Academy of Sciences, in
the summer of 2008 and during a workshop at the American Institute of Mathematics in Palo Alto in 
August of 2010.
They are indebted to both Institutes.}
\thanks{E. Werner has been partially supported by an NSF grant, a FRG-NSF grant and  a BSF grant.}}

\maketitle

\begin{abstract}
\noindent Let $K$ be a convex body in $\mathbb{R}^n$ with
Santal\'o point at $0$. We show that if  $K$ has a point
on the boundary with positive generalized Gau{\ss} curvature,
then the volume product $|K| |K^\circ|$ is not minimal.
This means that a body with minimal volume product has 
Gau{\ss} curvature equal to $0$ almost everywhere 
and thus suggests strongly that a minimal body
is a polytope.
\end{abstract}

\section{Introduction}

A convex body $K$ in $\mathbb R^{n}$ is a compact, convex set with
nonempty interior. The polar body $K^{z}$ with respect to an
interior point $z$ of $K$ is
$$
K^z= \{y |  \  \forall  x\in K: \  \langle y, x-z \rangle \leq 1\}.
$$
There is a unique point $z\in K$ such that the volume product
$|K| |K^\circ|$ is minimal.
This point is called the
Santal\'o point $s(K)$.
The Blaschke-Santal\'o inequality asserts that the maximum
of the volume product
$|K| |K^\circ|$ is attained for all ellipsoids and only for ellipsoids 
\cite{Bla,S,MP}. Thus the convex body for which
the maximum is attained is unique up to affine transforms.
\par
On the other hand, it is an open problem for which convex bodies the
minimum
is attained.  It is conjectured that the minimum is attained
for the simplex. The class of centrally symmetric convex bodies
is of particular importance. Mahler conjectured \cite{Ma1,Ma2} that
the minimum in this class is attained for the cube and its polar body, the
cross-polytope. If so,  the minimum would also be attained by
``mixtures'' of the cube and the cross-polytope,  sometimes called  Hanner-Lima bodies.
Those are not
affine images of the cube or the cross-polytope. Thus, in the class
of centrally symmetric convex bodies the minimum
is not attained for a unique convex body (up to affine transforms).
\par
The first breakthrough towards Mahler's conjecture is the inequality of Bourgain-Milman
\cite{BM}. They proved that for centrally symmetric convex bodies
$$
\left(\frac{c}{n}\right)^{n}
\leq
|K| |K^\circ|.
$$
This inequality has recently been reproved with completely
different methods by Kuperberg and Nazarov \cite{Ku,Naz}. Their proofs
also give better constants.
\par
For special classes like zonoids and unconditional
bodies Mahler's conjecture has been verified \cite{R1,GMR,S-R,M,R2}.
\par
The inequality of Bourgain-Milman has many applications in
various fields of mathematics: geometry of numbers, Banach space theory, 
convex geometry, theoretical computer science.
\par
Despite great efforts, a proof of Mahler's conjecture seems
still elusive. It is not even known whether a convex body
for which the minimum is attained must be a polytope.
A result in this direction has been proved by Stancu \cite{St}. It
is shown there that if $K$ is of class $C^2$ with strictly positive Gau{\ss} curvature everywhere,
then the volume product of $K$ can not be a local minimum.
\par
In this paper we show that a minimal body can not have even  a single  point with
positive generalized curvature. By a result of Alexandrov, Busemann and Feller \cite{Ale,BF}  the generalized curvature exists
almost everywhere. Therefore, our result implies that a minimal body has 
almost everywhere
curvature equal to $0$ and thus suggests strongly that a minimal body
is a polytope.

\vskip 4mm
We now introduce the concept of generalized curvature.
Let $\mathcal U$ be a convex, open subset of $\mathbb R^{n}$ and let
$f:\mathcal U \rightarrow \mathbb R$ be a convex function.
$d f(x) \in  \Bbb R^{n}$ is called
subdifferential at the point $x_{0} \in \mathcal U$, if we have for all
$x \in \mathcal U$
$$
f(x_{0})+\langle d f(x_{0}),x-x_{0}\rangle \leq f(x).
$$
A convex function has a subdifferential at every point and it is
differentiable at a point if and only if the subdifferential is unique.
Let $\mathcal U$ be an open, convex subset in $\mathbb R^{n}$ and $f:\mathcal
U\rightarrow\mathbb R$ a convex function. $f$ is said to be twice differentiable in
a  generalized sense in $x_{0}\in\mathcal U$, if there is a
linear map $d^{2}f(x_{0}):\mathbb R^n\to \mathbb R^n$
and a  neighborhood
$\mathcal{U}(x_{0})\subseteq\mathcal U$ such that we have for all $x \in
\mathcal{U}(x_{0})$ and for all subdifferentials $d f(x)$
$$
\|d f(x)-d f(x_{0})-d^{2}f(x_{0})(x-x_{0})\| \leq
\Theta(\|x-x_{0}\|)\|x-x_{0}\|.
$$
Here,  $\| \   \|$ is the standard Euclidean norm on $\mathbb{R}^n$ and $\Theta$ is a monotone function with $\lim_{t \to 0} \Theta(t)=0$.
$d^{2}f(x_{0})$ is called (generalized) Hesse-matrix.
If $f(0)=0$ and $d f(0)=0$ then we call the set
$$
\{x\in\mathbb R^{n}|x^{t}d^{2}f(0)x=1\}
$$
the indicatrix of Dupin at $0$. Since $f$ is convex,  this set is an ellipsoid
or a cylinder with a base that is an ellipsoid of lower dimension.
The eigenvalues of $d^{2}f(0)$ are
called generalized principal curvatures and their product is called the
generalized Gau{\ss}-Kronecker
curvature $\kappa$.  
\par
It  will always be this generalized Gau{\ss}  curvature that we mean throughout  the rest of the paper
though we may occasionally just call it Gau{\ss}  curvature.
Geometrically the eigenvalues of $d^{2}f(0)$
that are different from $0$ are
the lengths of the principal axes of the indicatrix raised to
the power $(-2)$.
\par
To define the generalized Gau{\ss} curvature $\kappa(x)$ of a convex body $K$ at a boundary point $x $   with  unique outer normal $N_K(x)$, 
if it exists, we translate and rotate $K$ so that we may assume that $x=0$ and $N_K(x)=-e_n$.
$\kappa(x)$ is then defined as the Gau{\ss} curvature of the function $f:\Bbb R^{n-1}\to \Bbb R$
whose graph in the neighborhood of $0$ is $\partial K$.

\par
We further denote by 
$H(x,\xi)$  the hyperplane through $x$ and orthogonal to $\xi$.  $H^-(x,\xi)$ and $H^+(x,\xi)$ are the two half spaces
determined by $H(x,\xi)$. In particular,
for $\Delta >0$, a convex body $K$ and $x\in\partial K$, the boundary of $K$,  with a unique
outer normal
$N_{K}(x)$
$$
H\left(x-\Delta N_K(x), N_K(x)\right)
$$ is the hyperplane through
$x-\Delta N_K(x)$ with  normal $N_K(x)$. $H^+\left(x-\Delta N_K(x), N_K(x)\right)$ denotes the halfspace determined by
$H\left(x-\Delta N_K(x), N_K(x)\right)$ that does not contain $x$.
\vskip 3mm
We construct two new bodies, $K_{x}(\Delta)$, by cutting off
a cap
\begin{equation*}\label{tilde}
K_{x}(\Delta)= K \cap H^+\left( x-\Delta N_K(x), N_K(x)\right),
\end{equation*}
and $K^{x}(\Delta)$ by
$$
K^{x}(\Delta)
=\operatorname{co} [K,x+\Delta N_{K}(x)].
$$

\vskip 4mm
\section{The main theorem}
\begin{thm} \label{mainthm}
Let $K$ be a convex body in $\mathbb{R}^n$ and suppose that there is a
point in the boundary of $K$ where the generalized Gau{\ss}  curvature
exists and is not $0$. Then the volume product
$|K|  | K^{s(K)}| $
is not
a local minimum. 
\par
Moreover, if $K$ is centrally symmetric with
center $0$ then, under the above assumption, the volume product $|K||K^{\circ}|$ is not
a local minimum in the class of $0$-symmetric convex bodies.
\end{thm}
\vskip 3mm
In order to prove Theorem \ref{mainthm},  we present the following
proposition.
\vskip 3mm

\begin{prop}\label{thm1}
Let $K$ be a convex body in $\mathbb{R}^n$ whose
Santal\'o point is at the origin. Suppose that there is a
point $x$ in the boundary of $K$ where the generalized Gau{\ss}
curvature exists and is not $0$. Then there exists $\Delta
>0$ such that
$$
|K_{x}(\Delta)| |\left(K_{x}(\Delta)\right)^\circ| < |K| |K^\circ|.
$$
or
$$
|K^{x}(\Delta)| |\left(K^{x}(\Delta)\right)^\circ| < |K| |K^\circ|.
$$
\end{prop}
\vskip 3mm

For the proof of Proposition \ref{thm1} we need several lemmas from \cite{SW4} and \cite{SW2004}.
We  refer to \cite{SW4} and \cite{SW2004} for the proofs. In particular, part (ii) of this lemma can be found in \cite{SW2004} as Lemma 12.

\vskip 3mm

\begin{lemma}\label{curvature}\cite{SW2004}
Let $K$ be a convex body in $\mathbb R^{n}$. Let $T:\mathbb R^{n}\rightarrow\mathbb R^{n}$ be a linear, invertible map.
\newline
(i)  The normal at $T(x)$ is
$$
(T^{-1})^t(N_{K}(x))\|(T^{-1})^t (N_{K}(x))\|^{-1}.
$$
(ii)
 Suppose that
the generalized Gau{\ss}-Kronecker curvature $\kappa$ exists in $x\in\partial K$.
Then the generalized Gau{\ss}-Kronecker curvature $\kappa$ exists in
$T(x)\in\partial T(K)$ and
$$
\kappa(x)^{}
=\|(T^{-1})^t(N_{K}(x))\|^{n+1}\det(T)^{2}
\kappa(T(x)).
$$
\end{lemma}
\vskip 4mm

The next two lemmas are well known. See e.g. \cite{SW4}.

\begin{lemma}\label{Lemma 1.2}\cite{SW4}
Let $\mathcal U$ be an open, convex subset of $\Bbb R^{n}$ and $0 \in \mathcal U$.
Suppose that $f:\mathcal U \rightarrow \Bbb R$ is twice differentiable in
the generalized sense at $0$ and that $f(0)=0$ and $d f(0)=0$.
\newline
Suppose that the indicatrix of Dupin at $0$ is an ellipsoid. Then there
is a monotone, increasing function $\psi:[0,1] \rightarrow [1,\infty )$ with
$\lim_{s \to 0}\psi(s)=1$ such that
\begin{eqnarray*}
& &\left\{(x,s)\left|x^{t}d^{2}f(0)x \leq \frac{2s}{\psi(s)}\right.\right\}
\\
& & \hskip 10mm\subseteq
\{(x,s)|f(x) \leq s\}
\subseteq
\{(x,s)|x^{t}d^{2}f(0)x \leq 2s \psi(s)\}.
\end{eqnarray*}
\end{lemma}
\vskip 4mm

\begin{lemma}\label{Lemma 1.3}\cite{SW4}
Let $K$ be a convex body in $\Bbb R^{n}$ with $0\in\partial K$ and
$N(0)=-e_{n}$. Suppose
that  the indicatrix of Dupin at $0$ is an ellipsoid. Suppose that the
principal axes $b_ie_i$ of the indicatrix are multiples of the unit
vectors $e_{i}$,
$i=1,\dots,n-1$. Let $\mathcal E$ be the $n$-dimensional ellipsoid
$$
\mathcal E=\left\{x\in\Bbb R^{n}\left|
\sum_{i=1}^{n-1}\frac{x_{i}^{2}}{b_{i}^{2}}
+\frac{\left(x_{n}-\left(\prod_{i=1}^{n-1}b_{i}\right)^{\frac{2}{n-1}}
\right)^{2}}
{(\prod_{i=1}^{n-1}b_{i})^{\frac{2}{n-1}}}\leq
\left(\prod_{i=1}^{n-1}b_{i}\right)^{\frac{2}{n-1}}\right.\right\}.
$$
Then there is an increasing, continuous function
$\phi:[0,\infty)\rightarrow [1,\infty)$ with $\phi(0)=1$ such that we have
for all $t$
\begin{eqnarray*}
& & \left\{\left.\left(\frac{x_{1}}{\phi(t)},\dots,\frac{x_{n-1}}{\phi(t)},
t\right)\right|x\in\mathcal
E, x_{n}=t\right\}   \\
& & \hskip 10mm\subseteq K\cap H((0,\dots,0,t),N(0))   \\
& & \hskip 10mm
\subseteq\left\{(\phi(t)x_{1},\dots,\phi(t)x_{n-1},t)|x\in\mathcal E,
x_{n}=t\right\}.
\end{eqnarray*}
We call $\mathcal E$ the standard approximating ellipsoid .
\end{lemma}
\vskip 4mm

Let us denote the lengths of the principal axes of the indicatrix of
Dupin by
$b_{i}$, $i=1,\dots,n-1$. Then the lengths $a_{i}$, $i=1,\dots,n$ of the principal
axes of the standard
approximating ellipsoid $\mathcal E$ are
\begin{equation}\label{1.3.1}
a_{i}=b_{i}\left(\prod_{i=1}^{n-1}b_{i}\right)^{\frac{1}{n-1}}
\hskip 5mm
i=1,\dots,n-1
\hskip 10mm
\mbox{and}
\hskip 5mm
a_{n}=\left(\prod_{i=1}^{n-1}b_{i}\right)^{\frac{2}{n-1}}.
\end{equation}
This follows immediately from Lemma \ref{Lemma 1.3}.
For the generalized Gau{\ss}-Kronecker curvature we get
\begin{equation}\label{1.3.2}
\prod_{i=1}^{n-1}\frac{a_{n}}{a_{i}^{2}}.
\end{equation}
This follows as the  generalized Gau{\ss}-Kronecker curvature equals the product of the
eigenvalues of the generalized Hesse matrix. The eigenvalues are $b_i^{-2}$, $i=1, \dots,
n-1$. Thus
$$
\prod_{i=1}^{n-1}b_{i}^{-2}
=\left(\prod_{i=1}^{n-1}b_{i}\right)^{2}
\prod_{i=1}^{n-1}\left(b_{i}\left(
\prod_{k=1}^{n-1}b_{k}\right)^{\frac{1}{n-1}}\right)^{-2}
=\prod_{i=1}^{n-1}\frac{a_{n}}{a_{i}^{2}}.
$$
In particular, if the indicatrix of Dupin is a sphere of radius
$\sqrt{\rho}$ then
the standard approximating ellipsoid is a Euclidean ball of radius $\rho$.
\par
We consider the map
$T:\Bbb R^{n}\rightarrow\Bbb R^{n}$
\begin{equation}\label{1.3.3}
T(x)=
\left(
\frac{x_{1}}{a_{1}}\left(\prod_{i=1}^{n-1}b_{i}\right)^{\frac{2}{n-1}},
\dots,\frac{x_{n-1}}{a_{n-1}}
\left(\prod_{i=1}^{n-1}b_{i}\right)^{\frac{2}{n-1}}, x_{n}\right).
\end{equation}
This transforms the standard approximating ellipsoid $\mathcal E$ into a
Euclidean ball $T(\mathcal E)$ with radius
$r=(\prod_{i=1}^{n-1}b_{i})^{2/(n-1)}$. This is
obvious since the principal axes of the standard approximating
ellipsoid are given by (\ref{1.3.1}). The map $T$ is volume preserving.
\vskip 4mm

\begin{lemma}\label{NormalPoint}
Let $K$ be a convex body in $\mathbb{R}^n$ with $0$
as an interior point. Suppose that the generalized
Gau{\ss} curvature of $\partial
K$ at $x_{0}$ exists and that the indicatrix is an
ellipsoid. Then there is an invertible  linear transformation
$T$ such that
\vskip 2mm
\noindent
(i) $N_{T(K)} (T(x_0))= \frac{T(x_0)}{\|T(x_0)\|}$
\vskip 2mm
\noindent
(ii) The indicatrix of Dupin at $T(x_0)$ is a Euclidean ball.
\vskip 2mm
\noindent
(iii)
$|T(K)|=|T(K)^\circ|$
\vskip 2mm
\noindent
(iv)
$\|T(x_0)\|=1$.
\end{lemma}

\vskip 3mm
\noindent
{\bf Proof.}
(i) We first show that there is a linear map $T_1$ such that $N_{T_1(K)} (T_1(x_0))= \frac{T_1(x_0)}{\|T_1(x_0)\|}$. Let $e_i$, $1 \leq i \leq n$ be an orthonormal basis of $\mathbb{R}^n$ such that $e_n=
N_K(x_0)$.
\par
\noindent
Define $T_{1}$ by
$$
T_{1}\left(\sum_{i=1}^{n-1}t_{i}e_{i}+t_{n}x_{0}\right)
=\sum_{i=1}^{n-1}t_{i}e_{i}+t_{n}\langle x_{0},e_{n}\rangle e_{n}
$$
$T_{1}$ is well-defined since $x_{0}$ and $e_{1},\dots,e_{n-1}$
are linearly independent. Indeed, as $\langle e_n, x_0 \rangle > 0$, $x_0 \notin e
^ \perp = \mbox{span}\{e_1, \dots, e_{n-1}\}$.
Moreover,
$$
T_{1}(x_{0})=\langle x_{0},e_{n}\rangle e_{n}
$$
and thus
\begin{equation}\label{en}
\frac{T_{1}(x_{0})}{\|T_{1}(x_{0}) \|} = e_{n}
\end{equation}
\noindent
By Lemma \ref{curvature},  the outer normal at $T_1(x_0)$ is
$$
(T_1^{-1})^t(N_{K}(x_{0}))\|(T_1^{-1})^t (N_{K}(x_0))\|^{-1}.
$$
Then for all $1 \leq i \leq n-1$
$$
\langle (T_1^{-1})^t(N_{K}(x_{0})), e_{i}\rangle
=\langle N_{K}(x_{0}), T_1^{-1}e_{i}\rangle
=\langle N_{K}(x_{0}),e_{i}\rangle
=\langle e_{n},e_{i}\rangle
=0.
$$
Hence
$$
\frac{(T_1^{-1})^t (N_{K}(x_{0}))}{\|(T_1^{-1})^t (N_{K}(x_0))\|} = \pm e_n
$$
Since $\langle x_{0},e_{n}\rangle >0$ and since $0$ is an interior point, it is $+ e_n$.
Together with (\ref{en}), this shows that
$$
N_{T_1(K)} (T_1(x_0)) = e_{n} = \frac{T_{1}(x_{0})}{\|T_{1}(x_{0}) \|}.
$$

\vskip 4mm
\noindent
(ii) Put $x_1=T_1(x_0)$ and $K_1=T_1(K)$. By Lemma \ref{curvature}, the curvature
$\kappa(x_1)$
at $x_1 \in \partial K_1$ exists and is positive.
For $1 \leq i \leq n-1$, let $b_{i}$ be the principal curvatures and $a_{i}$ be the principal axes of the standard approximating ellipsoid in $x_1$.
Let
$T_2:\Bbb R^{n}\rightarrow\Bbb R^{n}$
\begin{equation}\label{1.3.3}
T_2(x)=
\left(
\frac{\xi_{1}}{a_{1}}\left(\prod_{i=1}^{n-1}b_{i}\right)^{\frac{2}{n-1}},
\dots,\frac{\xi_{n-1}}{a_{n-1}}
\left(\prod_{i=1}^{n-1}b_{i}\right)^{\frac{2}{n-1}},\xi_{n}\right).
\end{equation}
transforms the indicatrix of Dupin at $x_1$ into an $n-1$-dimensional Euclidean ball and
the standard approximating ellipsoid $\mathcal E$ into a $n$-dimensional
Euclidean ball $T_2(\mathcal E)$ with radius
$r=(\prod_{i=1}^{n-1}b_{i})^{2/(n-1)}$.
\par
Property (i) of the lemma is preserved:
$$
N_{T_{2}(K_{1})}(T_{2}(x_{1})) = N_{K_{1}}(x_{1})=e_{n}.
$$
Indeed, by Lemma \ref{curvature}
$$
N_{T_{2}(K_{1})}(T_{2}(x_{1}))
=\frac{(T_{2}^{-1})^t (N_{K_{1}}(x_{1}))}{\|(T_{2}^{-1})^t(N_{ K_{1}}(x_{1}))\|}.
$$
and thus for all $1 \leq i \leq n-1$
$$
\langle (T_{2}^{-1})^t (N_{K_{1}}(x_{1})),e_{i}\rangle
=\langle N_{K_{1}}(x_{1}), T_{2}^{-1}(e_{i})\rangle
=\left\langle e_{n},a_{1}\left(\prod_{i=1}^{n-1}b_{i}\right)^{-\frac{2}{n-1}}e_{i}\right\rangle
=0.
$$
\vskip 4mm
\noindent
(iii)
It is enough to apply a multiple $\alpha I$ of the identity.
\vskip 4mm
\noindent
(iv) We apply the map $T_{3}$ with
$$
T_{3}(\xi)
=(\lambda \xi_{1},\dots,\lambda \xi_{n-1},\lambda^{-n+1}\xi_{n})
$$
where
$$
\lambda
= (\alpha \langle x_0, e_n \rangle )^{1/(n-1)}
$$
Properties (i) and (ii)  of the lemma are preserved and, as $\mbox{det}(T_{3})=1$, Property (iii) as well.
\par
\noindent
Finally,  we let $T(K)=T_3(\alpha T_2(K_1))$.
$\Box$

\vskip 4mm

\begin{lemma}\label{DiffDual}
Let $K$ be a convex body in $\mathbb R^{n}$ such that
$\partial K$ is twice differentiable in the generalized
sense at $x$. Suppose that $\|x\|=1$, $N_{ K}(x)=x$,
and the indicatrix of Dupin at $x$   is a Euclidean sphere with
radius
$r$.
 Then $x\in\partial K^{\circ}$
and for all $0 < \epsilon < \operatorname{min}\{ r, \frac{1}{r}\}$ there is $\Delta>0$ such that
\begin{eqnarray*}
&& \hskip -10mm B_{2}^{n}\left(x-\left(\frac{1}{r}-\epsilon\right)N_{
K^{\circ}}(x), \frac{1}{r}-\epsilon \right)\cap
H^{-}(x-\Delta N_{
K^{\circ}}(x),N_{
K^{\circ}}(x))    \\
&& \subseteq K^{\circ} \cap
H^{-}(x-\Delta N_{
K^{\circ}}(x),N_{
K^{\circ}}(x)) \\
&&\hskip 10mm  \subseteq
B_{2}^{n}\left(x-\left(\frac{1}{r}+\epsilon\right)N_{
K^{\circ}}(x), \frac{1}{r}+\epsilon \right)\cap
H^{-}(x-\Delta N_{
K^{\circ}}(x),N_{
K^{\circ}}(x))
\end{eqnarray*}
\end{lemma}
\vskip 3mm

\noindent
{\bf Proof.}
Without loss of generality we can assume that $x =N_{K(x)}= e_n$.
Clearly then $x \in \partial K^{\circ}$ and $N_{K^{\circ}(x)}=x$.
Let $0 < \epsilon < \operatorname{min}\{ r, \frac{1}{r}\}$.  By Lemma \ref{Lemma 1.3}, there exists $\Delta_1$ such that for all $\Delta \leq \Delta_1$
\begin{eqnarray*}
&&B_{2}^{n}\left((1-(r-\epsilon)) e_n, r-\epsilon \right)
\cap
H^{-}\left((1-\Delta) e_n, e_n\right)    \\
&&\subseteq K \cap H^{-}\left((1-\Delta) e_n, e_n\right)
\subseteq
B_{2}^{n}\left((1-(r+\epsilon)) e_n, r+\epsilon \right) \cap
H^{-}\left((1-\Delta) e_n, e_n\right).
\end{eqnarray*}
We construct  two new convex  bodies.
\begin{eqnarray*}
K_1= \operatorname{co}\left[ K \cap H^{+}\left((1-\Delta_1) e_n, e_n\right),  B_{2}^{n}\left((1-(r-\epsilon)) e_n, r-\epsilon \right)\right]
\end{eqnarray*}
and
\begin{eqnarray*}
K_2= \operatorname{co}\left[ K \cap H^{+}\left((1-\Delta_1) e_n, e_n\right),  B_{2}^{n}\left((1-(r+\epsilon)) e_n, r+\epsilon \right)\right].
\end{eqnarray*}
Then $K_1 \subseteq K \subseteq K_2$ and there is $\Delta_2 \leq \Delta_1$ such that for all $\Delta \leq \Delta_2$
\begin{eqnarray*}
K_1 \cap H^{-}\left((1-\Delta) e_n, e_n\right)= B_{2}^{n}\left((1-(r-\epsilon)) e_n, r-\epsilon \right) \cap H^{-}\left((1-\Delta) e_n, e_n\right)
\end{eqnarray*}
and
\begin{eqnarray*}
K_2 \cap H^{-}\left((1-\Delta) e_n, e_n\right)= B_{2}^{n}\left((1-(r+\epsilon)) e_n, r-\epsilon \right) \cap H^{-}\left((1-\Delta) e_n, e_n\right).
\end{eqnarray*}
We now compute $K_1^0$  and $K_2^0$ in a neighborhood of $x=e_n$. We show the computations for $K_1$. $K_2$ is done
similarly.
\par
Let  $\Delta  \leq \Delta_2$ and $\eta$ be the normal
of  $y \in \partial K_1 \cap H^{-}\left((1-\Delta) e_n, e_n\right)$. Then
$$
\langle y,\eta \rangle
=r- \epsilon +\left(1-(r-\epsilon)\right) \  \langle x,\eta \rangle
$$
Therefore,
$$
\bigg\langle  y,\frac{\eta}{r- \epsilon+(1-(r-\epsilon) ) \  \langle x,\eta \rangle } \bigg\rangle=1
$$
and hence
$$
\frac{\eta}{r-\epsilon+(1-(r-\epsilon) ) \langle x,\eta \rangle }  \in \partial K_1^{\circ}.
$$
For  $\Delta_{K_1^\circ} \leq \Delta_2$ sufficiently small, we consider now a cap of $K_1^{\circ}$ and its base
$$
K_1^{\circ}
\cap H((1-\Delta_{K_1^{\circ}})e_n, e_n)
$$
We compute the distance $\rho$ of
$\frac{\eta}{r - \epsilon+(1-(r-\epsilon) ) \langle x,\eta \rangle} \in  \partial K_1^{\circ}
\cap H((1-\Delta_{K_1^{\circ}})e_n, e_n)$ from the center of the
base $(1-\Delta_{K_1^{\circ}})e_n$ . Clearly,
by Pythagoras
$$
\rho
=\sqrt{\frac{1}{(r-\epsilon+(1-(r-\epsilon) ) \langle x,\eta \rangle)^{2}}
-(1-\Delta_{K_1^{\circ}})^{2}}
$$
Moreover,
$$
1-\Delta_{K_1^{\circ}}
=\frac{\langle x,\eta \rangle}{r-\epsilon+[1-(r -\epsilon)]\langle x,\eta\rangle}
$$
Hence
$$
(1-\Delta_{K_1^{\circ}}) \  \bigg(r -\epsilon +[1-(r-\epsilon)] \ \langle x,\eta\rangle \bigg)
= \langle x,\eta\rangle
$$
$$
\langle x,\eta\rangle  \bigg((1-\Delta_{K_1^{\circ}}) [1-(r-\epsilon)]-1\bigg)
= - (r-\epsilon) (1-\Delta_{K_1^{\circ}})
$$
$$
\langle x,\eta \rangle
=\frac{(r-\epsilon) (1-\Delta_{K_1^{\circ}})}
{r-\epsilon+[1-(r-\epsilon)] \ \Delta_{K_1^{\circ}}}
$$
Therefore
$$
\rho
=\sqrt{\frac{\left(r-\epsilon+[1-(r-\epsilon)]\Delta_{K_1^{\circ}}\right)^{2}}
{(r-\epsilon)^{2}}
-(1-\Delta_{K_1^{\circ}})^{2}}
$$
or
$$
\rho
=\sqrt{\frac{2}{r-\epsilon}\Delta_{K_1^{\circ}}
+\frac{[1-2(r-\epsilon)]}{(r-\epsilon)^{2}}\  \Delta_{K_1^{\circ}}^{2}}
$$
We compare this radius with the corresponding radius of the ball
$B_{2}^{n}\left((1-\frac{1}{r-2\epsilon})e_{n},\frac{1}{r-2\epsilon}\right)$.
The corresponding radius is
$$
\sqrt{\frac{2\Delta_{K_1^{\circ}}}{r-2\epsilon}-\Delta_{K_1^{\circ}}^{2}}
$$
Therefore,
$$
\rho\leq\sqrt{\frac{2\Delta_{K_1^{\circ}}}{r-2\epsilon}-\Delta_{K_1^{\circ}}^{2}}
$$
provided that
$$
\Delta_{K_1^{\circ}}\leq\min\left\{\frac{2\epsilon}{(r-2\epsilon)(r-\epsilon)},
\frac{2\epsilon}{(r-2\epsilon)(2+r-\epsilon)}\right\}.
$$
As $K_{1}\subseteq K$ and consequently $K^{\circ}\subseteq K_{1}^{\circ}$
we get
$$
K^{\circ}\cap H^{-}((1-\Delta_{K_{1}^{\circ}})e_{n},e_{n})
\subseteq
B_{2}^{n}\left(\left(1-\frac{1}{r-2\epsilon}\right)e_{n},\frac{1}{r-2\epsilon}\right)
\cap H^{-}((1-\Delta_{K_{1}^{\circ}})e_{n},e_{n}).
$$
Similarly, using  $K \subseteq K_2$, one shows (with a new  $\Delta _{K_{1}^{\circ}}$ small enough if needed) that
$$
B_{2}^{n}\left(\left(1-\frac{1}{r+2\epsilon}\right)e_{n},\frac{1}{r+2\epsilon}\right)
\cap H^{-}((1-\Delta_{K_{1}^{\circ}})e_{n},e_{n})  \subseteq
K^{\circ}\cap H^{-}((1-\Delta_{K_{1}^{\circ}})e_{n},e_{n}).
$$
$\Box$
\vskip 4mm
Lemma \ref{DiffDual} implies that if $x \in \partial K$  is a twice differentiable point (in the generalized sense), then the point $y \in \partial K^\circ$ with $\langle x, y \rangle  =1$, is also a twice differentiable point.
Compare also  \cite{Hug}.
\begin{cor} 
Let $K$ be a convex body in $\mathbb R^{n}$ with $0$ as an interior
point. Assume that
$\partial K$ is twice differentiable in the generalized
sense at $x$ and the indicatrix of Dupin at $x$ is an ellipsoid
(and not a cylinder with an ellipsoid as its base). Then $\partial K^{\circ}$
is twice differentiable at the unique point $\xi$ with 
$\langle \xi,x\rangle=1$.
\end{cor}
\vskip 3mm

\noindent
{\bf Proof.}
By Lemma \ref{NormalPoint} we may assume that the indicatrix
is a Euclidean ball and that $N_{K}(x)=x$. By Lemma \ref{DiffDual}
the statement follows. 
$\Box$
\vskip 4mm

The next lemma is also well known (see e.g. \cite{SW1990}).
\begin{lemma}\label{cap}
Let $K$ be a convex body in $\mathbb R^{n}$
and suppose that the indicatrix of Dupin at
$x\in\partial K$ exists and is a Euclidean ball of
radius $r>0$. Let $C(r,\Delta)$ be the cap at $x$
of height $\Delta$. Then
$$
|C(r,\Delta)|
=g\left(\frac{\Delta}{r}\right)^{\frac{n+1}{2}}
\frac{1}{n+1}2^{\frac{n+1}{2}}
\operatorname{vol}_{n-1}(B_{2}^{n-1})
\Delta^{\frac{n+1}{2}}r^{\frac{n-1}{2}}
$$
where $\lim_{t\to0}g(t)=1$.
\end{lemma}
\vskip 4mm
\noindent
{\bf Remark}. The conclusions of Lemma \ref{cap} also hold if instead of the existence of the indicatrix, we assume
the following: 
\newline
Let  $x \in \partial K$ 
and suppose that there is $r >0$ such that
for all $\epsilon>0$ there is a
$\Delta_{\epsilon}$ such that for all $\Delta$ with
$0<\Delta\leq\Delta_{\epsilon}$
\begin{eqnarray} \label{approx}
&& \hskip -10mm B_{2}^{n}\left(x-\left( r-\epsilon\right)N_{
K}(x), r-\epsilon \right)\cap
H^{-}(x-\Delta N_{
K}(x),N_{
K}(x))  \nonumber  \\
&& \subseteq K\cap
H^{-}(x-\Delta N_{
K}(x),N_{
K}(x))\nonumber \\
&&\hskip 10mm  \subseteq
B_{2}^{n}\left(x-\left(r+\epsilon\right)N_{
K}(x), r+\epsilon \right)\cap
H^{-}(x-\Delta N_{
K}(x),N_{
K}(x)).
\end{eqnarray}

\vskip 4mm
\noindent
The next lemma is from  \cite{W1}.  There  it was  assumed that  the indicatrix of Dupin at
$x\in\partial K$ exists and is a Euclidean ball of
radius $r>0$. However, what was actually used in the proof, were the assumptions (\ref{approx}) of the above Remark.
\par
\begin{lemma}\label{cone1}\cite{W1}
Let $K$ be a convex body in $\mathbb R^{n}$. Let $x \in \partial K$ 
and suppose that there is $r >0$ such that
for all $\epsilon>0$ there is a
$\Delta_{\epsilon}$ such that for all $\Delta$ with
$0<\Delta\leq\Delta_{\epsilon}$, (\ref{approx}) holds.
Then,  if $\Delta_{\epsilon}$ is small enough, we have
for $0<\Delta<\Delta_{\epsilon}$
\begin{eqnarray*}
&&\hskip -10mm\frac{2^{\frac{n+1}{2}}}{n(n+1)}
|B_{2}^{n-1}|\Delta^{\frac{n+1}{2}}
r^{\frac{n-1}{2}} \\
&& \times \left\{(n+1)
\left(1-\frac{\epsilon}{r}\right)^{\frac{n-1}{2}}
(1-c\  \epsilon)
-n\left(1+\frac{\epsilon}{r}\right)^{\frac{n-1}{2}}
(1+c\ \epsilon)\  h \left(\frac{\Delta}{r+\epsilon}
\right)^{n+1}
\right\}   \\
&& \hskip -10mm \leq
 \bigg| K^{x}(\Delta)\setminus K \bigg|
\\
&&\hskip -10mm\leq\frac{2^{\frac{n+1}{2}}}{n(n+1)}
|B_{2}^{n-1}|\Delta^{\frac{n+1}{2}}
r^{\frac{n-1}{2}} \\
&& \times \left\{(n+1)
\left(1+\frac{\epsilon}{r}\right)^{\frac{n-1}{2}}
(1+ c \ \epsilon)
-n\left(1-\frac{\epsilon}{r}\right)^{\frac{n-1}{2}}
(1-c\  \epsilon) \  h\left(\frac{\Delta}{r+\epsilon}
\right)^{n+1}
\right\}.
\end{eqnarray*}
where $c$ is a constant and $\lim_{t\to0}h(t)=1$.
\end{lemma}
\vskip 5mm

\noindent
{\bf Proof of Theorem \ref{mainthm}.}
\vskip 3mm
\noindent
By assumption there is a point $x\in\partial K$
at which $\partial K$ is twice differentiable in the
generalized sense. By Lemma \ref{NormalPoint}
we may assume that $\|x\|=1$ and $x=N_{K}(x)$.
Moreover, all principal curvature radii at $x$ are
equal to $r$.
By Lemma \ref{DiffDual},  $x\in\partial K^{\circ}$,
$K^{\circ}$ is twice differentiable at $x$
and all principal curvature radii are equal to
$\frac{1}{r}$.
\par
The dual body to $K_{x}(\Delta)$ is
$K^{x}(\frac{\Delta}{1-\Delta})$.
By Lemma \ref{cap},
$$
|K_{x}(\Delta)|
\leq|K|
- g\left(\frac{\Delta}{r}\right)^{\frac{n+1}{2}}  \frac{2^{\frac{n+1}{2}}}{n+1}
\operatorname{vol}_{n-1}(B_{2}^{n-1})
\Delta^{\frac{n+1}{2}}r^{\frac{n-1}{2}}.
$$
By Lemma \ref{cone1}
\begin{eqnarray*}
&& \hskip -15mm | K_{x}(\Delta)^{\circ} |
= \left|K^{x}\left(\frac{\Delta}{1-\Delta}\right)\right|
\leq |K^{\circ}|
+ \frac{2^{\frac{n+1}{2}}}{n(n+1)}
|B_{2}^{n-1}|\Delta^{\frac{n+1}{2}}
r^{-\frac{n-1}{2}} \\
&& \times \left\{(n+1)
\left(1+r\epsilon\right)^{\frac{n-1}{2}}
(1+ c\  \epsilon)
-n\left(1-r\epsilon\right)^{\frac{n-1}{2}}
(1- c\ \epsilon)\  h\left(\frac{\Delta}{\frac{1}{r}+\epsilon}
\right)^{n+1}
\right\}.
\end{eqnarray*}
It follows that
\begin{eqnarray*}
&& \hskip -15mm |K_{x}(\Delta)|
|K_{x}(\Delta)^{\circ}|
\leq|K||K^{\circ}|
+|K|\frac{2^{\frac{n+1}{2}}}{n(n+1)}
|B_{2}^{n-1}|\Delta^{\frac{n+1}{2}}
r^{-\frac{n-1}{2}} \\
&& \times \left\{(n+1)
\left(1+r\epsilon\right)^{\frac{n-1}{2}}
(1+ c\ \epsilon)
-n\left(1-r\epsilon\right)^{\frac{n-1}{2}}
(1-c\ \epsilon)\  h\left(\frac{\Delta}{\frac{1}{r}+\epsilon}
\right)^{n+1}
\right\}   \\
&&-|K^{\circ}|g\left(\frac{\Delta}{r}\right)^{\frac{n+1}{2}}
\frac{1}{n+1}2^{\frac{n+1}{2}}
\operatorname{vol}_{n-1}(B_{2}^{n-1})
\Delta^{\frac{n+1}{2}}(r-\epsilon)^{\frac{n-1}{2}}.
\end{eqnarray*}
Therefore we have
$$
|K_{x}(\Delta)|
|K_{x}(\Delta)^{\circ}|
<|K||K^{\circ}|,
$$
provided that
\begin{eqnarray}\label{ineq1}
&&|K|\left\{(n+1)
\left(1+r\epsilon\right)^{\frac{n-1}{2}}
(1+c\ \epsilon)
-n\left(1-r\epsilon\right)^{\frac{n-1}{2}}
(1-c\ \epsilon)\  h \left(\frac{\Delta}{\frac{1}{r}+\epsilon}
\right)^{n+1}
\right\}   \nonumber \\
&&< n|K^{\circ}| \  g\left(\frac{\Delta}{r}\right)^{\frac{n+1}{2}}r^{n-1}
\left(1-\frac{\epsilon}{r}\right)^{\frac{n-1}{2}}.
\end{eqnarray}
Now we interchange the roles of $K$ and $K^{\circ}$.
We cut off a cap from $K^{\circ}$ and  apply the remark following
Lemma \ref{cap}.
Then
the inequality analogous to (\ref{ineq1}) will be
\begin{eqnarray} \label{ineq2}
&&|K^{\circ}|\left\{(n+1)
\left(1+\frac{\epsilon}{r}\right)^{\frac{n-1}{2}}
(1+c\ \epsilon)
-n\left(1-\frac{\epsilon}{r}\right)^{\frac{n-1}{2}}
(1-c\ \epsilon)\  h \left(\frac{\Delta}{r+\epsilon}
\right)^{n+1}
\right\} \nonumber   \\
&&<n|K|g\left(r\Delta\right)^{\frac{n+1}{2}}r^{-(n-1)}
(1-r\epsilon)^{\frac{n-1}{2}}.
\end{eqnarray}
Thus the theorem is proved provided
that one of the inequalities (\ref{ineq1}) or (\ref{ineq2}) holds.
Suppose  both inequalities do not hold. Then
\begin{eqnarray*}
&&|K|\left\{(n+1)
\left(1+r\epsilon\right)^{\frac{n-1}{2}}
(1+c\ \epsilon)
-n\left(1-r\epsilon\right)^{\frac{n-1}{2}}
(1-c\ \epsilon)\  h \left(\frac{\Delta}{\frac{1}{r}+\epsilon}
\right)^{n+1}
\right\}    \\
&&\geq n|K^{\circ}| \  g\left(\frac{\Delta}{r}\right)^{\frac{n+1}{2}}r^{n-1}
\left(1-\frac{\epsilon}{r}\right)^{\frac{n-1}{2}}   \\
&&\geq n^{2}|K|
\frac{g\left(\frac{\Delta}{r}\right)^{\frac{n+1}{2}}
(1-\frac{\epsilon}{r}-r\epsilon+\epsilon^{2})^{\frac{n-1}{2}}
g\left(r\Delta\right)^{\frac{n+1}{2}}
}{\left\{(n+1)
\left(1+\frac{\epsilon}{r}\right)^{\frac{n-1}{2}}
(1+c\ \epsilon)
-n\left(1-\frac{\epsilon}{r}\right)^{\frac{n-1}{2}}
(1-c\ \epsilon)\  h\left(\frac{\Delta}{r+\epsilon}
\right)^{n+1}
\right\}  }.
\end{eqnarray*}
We can choose $\epsilon$ so small that
$$
(n+1)
\left(1+r\epsilon\right)^{\frac{n-1}{2}}
(1+c\ \epsilon)
-n\left(1-r\epsilon\right)^{\frac{n-1}{2}}
(1-c\ \epsilon)\ \  h\left(\frac{\Delta}{\frac{1}{r}+\epsilon}
\right)^{n+1}\leq 2
$$
and
$$
(n+1)
\left(1+\frac{\epsilon}{r}\right)^{\frac{n-1}{2}}
(1+c\  \epsilon)
-n\left(1-\frac{\epsilon}{r}\right)^{\frac{n-1}{2}}
(1-c\ \epsilon)\ \  h\left(\frac{\Delta}{r+\epsilon}
\right)^{n+1}
  \leq2.
$$
Moreover, we can choose $\epsilon$ so small that
$$
\left(1-\frac{\epsilon}{r}-r\epsilon+\epsilon^{2}
\right)^{\frac{n-1}{2}}
\geq\frac{1}{2}
$$
Therefore
$$
4\geq n^{2}g\left(\frac{\Delta}{r}\right)^{\frac{n+1}{2}}
g\left(r\Delta\right)^{\frac{n+1}{2}}
$$
Since $\lim_{t\to0}g(t)=1$,  this gives a contradiction.
\par
The extension of the proof  needed in order to prove the symmetric case is obvious.
$\Box$

\bigskip
\noindent Shlomo Reisner   \\
{\small Department of Mathematics}\\
{\small University of Haifa }\\
{\small Haifa, 31905 Israel}\\
{\small \tt reisner@math.haifa.ac.il}\\ \\
\vskip 2mm
\and
\noindent
Carsten Sch\"utt  \\
{\small Mathematisches Seminar }\\
{\small Christian-Albrechts Universit\"at }\\
{\small D-24098 Kiel, Ludewig-Meyn Strasse 4, Germany}\\
{\small \tt schuett@math.uni-kiel.de}\\ \\

\and
\vskip 2mm
\noindent Elisabeth Werner\\
{\small Department of Mathematics \ \ \ \ \ \ \ \ \ \ \ \ \ \ \ \ \ \ \ Universit\'{e} de Lille 1}\\
{\small Case Western Reserve University \ \ \ \ \ \ \ \ \ \ \ \ \ UFR de Math\'{e}matique }\\
{\small Cleveland, Ohio 44106, U. S. A. \ \ \ \ \ \ \ \ \ \ \ \ \ \ \ 59655 Villeneuve d'Ascq, France}\\
{\small \tt elisabeth.werner@case.edu}\\ \\


\vskip 3mm



\begin{thebibliography}{00}

\bibitem{Ale}
A.D. Alexandrov, 
{\it Almost everywhere existence of the second differential of a convex
function and some properties of convex surfaces connected with it},
Uchen. Zap. Leningrad Gos. Univ. Mat. Ser. {\bf 6} (1939), 3-35.

\bibitem{Bla}
W. Blaschke, {\it \"Uber affine Geometrie VII. Neue Extremeigenschaften von Ellipse
und Ellipsoid\/}, Leipz. Ber., {\bf 69} (1917), 306-318.

\bibitem{BM}
J. Bourgain and  V. D. Milman, {\it New volume ratio properties for
convex symmetric bodies in $\Bbb R^n$\/}, Invent.\ Math.\ {\bf 88} (1987),
319-340.

\bibitem{BF}
H. Busemann and W. Feller,
{\it Kr\"ummungseigenschaften konvexer Fl\"achen},
Acta Math. {\bf 66} (1936), 1-47.


\bibitem{GMR}
Y. Gordon, M. Meyer and S. Reisner, {\it Zonoids with minimal
volume product--a new proof\/}, Proc.\ Amer.\ Math.\ Soc.\ {\bf 104} (1988),
273-276.

\bibitem{Hug}
{\sc D. Hug}, {\em Curvature Relations and Affine Surface Area for
a General Convex Body and its Polar},  Results in Mathematics  {\bf
V. 29}  (1996), 233-248.


\bibitem{Ku}
G. Kuperberg, {\it  From the Mahler Conjecture to Gau{\ss}
Linking Integrals\/}, Geometric And Functional Analysis {\bf 18} (2008), 870-892.

\bibitem{Ma1}
K. Mahler, {\it Ein Minimalproblem f\"ur konvexe Polygone\/},
Mathematica (Zutphen) {\bf B 7} (1939), 118-127.

\bibitem{Ma2}
K. Mahler, {\it Ein \"Ubertragungsprinzip f\"ur konvexe K\"orper\/},
 \v{C}asopis P\v{e}st.\ Mat.\ Fys.\ {\bf 68} (1939), 93-102.

\bibitem{M}
M. Meyer, {\it Une caract\'erisation volumique de certains
espaces norm\'es\/}, Israel J. Math.\ {\bf 55} (1986), 317-326.

\bibitem{MP}
M. Meyer and A. Pajor, {\it On the Blaschke-Santal\'o inequality\/}, 
Arch.\ Math.\ {\bf 55} (1990), 82-93.

\bibitem{Naz}
F. Nazarov, {\it The H\"ormander proof of the Bourgain-Milman theorem\/}, preprint, 2009.

\bibitem{R1}
S. Reisner, {\it Zonoids with minimal volume-product\/}, Math.\ Z.\ {\bf 192} (1986),
339-346.

\bibitem{R2}
S. Reisner,  {\it Minimal volume product in Banach spaces with
a 1-unconditional basis\/}, J. London Math.\ Soc.\ {\bf 36} (1987), 126-136.

\bibitem{S-R}
J. Saint-Raymond, {\it Sur le volume des corps convexes
sym\'etriques\/},  S\'eminaire d'Initiation \`a l'Analyse, 1980-1981, Universit\'e
PARIS VI, Paris 1981.

\bibitem{S}
L.\ A.\ Santal\'o, {\it Un invariante afin para los cuerpos convexos del
espacio de $n$ dimensiones\/}, Portugal.\ Math.\ {\bf 8} (1949), 155-161.


\bibitem{SW1990}
{C. Sch\"{u}tt and  E. Werner}, {\em The convex floating body\/},
Math. Scand. {\bf 66} (1990), 275-290.


\bibitem{SW4}
{C. Sch{\"u}tt and E. Werner}, {\em Random polytopes of points
chosen from the boundary of a convex body\/} GAFA Seminar
Notes, Lecture Notes in Math., {\bf 1807}, Springer-Verlag,
(2002),  241-422.

\bibitem{SW2004}
{C. Sch\"utt and E. Werner},
{\em Surface bodies and p-affine surface area},
Advances in Mathematics {\bf 187}  (2004), 98-145.

\bibitem{St}
A. Stancu, {\it Two volume product inequalities and their applications\/},
Canad. Math. Bull. {\bf 52} (2009), 464-472.


\bibitem{W1}
{ E. Werner},
{\em Illumination bodies and affine surface area},
Studia Mathematica {\bf 110} (1994), 257-269.


\end{thebibliography}
\end{document}